\documentclass[11pt,a4paper,reqno]{amsart}
\usepackage{mathrsfs,amssymb}

\newcommand{\N}{\mathbb{N}}

\newtheorem{theorem}{Theorem}

\newtheorem{corollary}{Corollary}

\newtheorem{proposition}{Proposition}
\newtheorem{lemma}{Lemma}

\newcounter{fiddletheoremtemp}

\begin{document}

\title[Generic Monoids and Semigroups]{On Generic Properties of Finitely Presented Monoids and Semigroups}

\maketitle

\begin{center}

    MARK KAMBITES

    \medskip

    School of Mathematics, \ University of Manchester, \\
    Manchester M13 9PL, \ England.

\end{center}

\begin{abstract}
We study the generic properties of finitely presented monoids and
semigroups. We show that for positive integers $a > 1$, $k$ and $m$,
the generic $a$-generator $k$-relation monoid and semigroup (defined in any of several
definite statistical senses) satisfy the small overlap
condition $C(m)$. It follows that the generic monoid is torsion-free and
$\mathscr{J}$-trivial and, by a recent result of the author, admits a linear
time solution to its word problem and a regular language of unique normal
forms for its elements. Moreover, the uniform word problem for finitely
presented monoids is generically solvable in time linear in the word lengths
and quadratic in the presentation size. We also prove some technical results
about generic sets which may be of independent interest.
\end{abstract}

\section{Introduction}

Traditional complexity theory studies the time taken to solve a problem or
execute an algorithm in the ``worst case'', but for many problems the
``worst case'' arises very infrequently. Probably the best known example
is Dantzig's simplex method for linear programming \cite{Dantzig47}, which has
exponential worst case time complexity but in practice almost invariably
terminates in linear time (see eg.~\cite{Klee69}). Now over 60 years old, it remains the preferred
choice for practical applications, even though there are now alternative algorithms with
worst-case polynomial time complexity. Phenomena such as this motivated
the development of
\textit{average-case complexity} \cite{Gurevich91}, which measures, roughly
speaking, the \textit{mean} difficulty of a problem across instances,
with respect to some measure. Average-case complexity has proved extremely
helpful for obtaining a theoretical understanding of the ``practical'' difficulty of problems,
especially within the class NP of problems admitting non-deterministic
worst-case polynomial time solution.

Average-case analysis can also be applied outside NP, but it
meets with a conceptual difficulty. For most applications, what matters
is not so much the \textit{mean} difficulty of a problem of across
instances, but rather the \textit{typical} difficulty of instances encountered
in practice. As is well-known to statisticians, the mean value of a data set
is not necessarily a guide to the typical values, since the former can be
heavily skewed in one direction by a very small number of very extreme
outliers. Likewise, the average-case complexity of a problem can be skewed
upwards by a very small proportion of very difficult instances. Within NP
worst cases are single exponential; this imposes a limit on the
``extremeness'' of outlying instances and hence their ability to distort
the mean. Outside NP, however, the distortion can be much more dramatic,
with a tiny minority of extremely difficult instances potentially inflating
the average-case complexity well beyond the complexity of the typical instance
encountered in practice. This culminates in the extreme case of recursively
unsolvable problems, whose average-case complexity is not defined at all,
even though algorithms may exist to solve such problems efficiently for an
overwhelming majority of cases \cite{Gilman07}.

The aim of \textit{generic-case complexity} is directly to analyse the
complexity of \textit{typical} problem instances, as distinct from the
\textit{average} difficulty of problem instances. Rather than introducing
a \textit{measure} on the instance space, the key idea is the
\textit{stratification} of an instance space (or indeed any other set)
into an infinite sequence of finite subsets. A subset $X$ of the space is
called \textit{generic} if the proportion of elements in each finite
set which belong to $X$ approaches $1$ as one moves along the sequence.
The generic complexity is (very roughly speaking) the minimum complexity
attainable on a generic set. Compared with the average-case approach, the
key feature is that no single instance (indeed no finite set of instances),
makes any contribution at all to the generic properties of the space. 
Generic-case complexity was introduced by group theorists \cite{Kapovich03},
investigating the large stock of hard algorithmic problems which occur in
the study of finitely generated infinite groups. It has proved especially
useful in view of recent interest in the use of non-commutative groups
as a basis for cryptographic systems \cite{Shpilrain06},
permitting for example a theoretical understanding of the success of the
\textit{length-based attack} \cite{Ruinskiy07} on the Shpilrain-Ushakov
key establishment protocol based on the Thompson group \cite{Shpilrain05}.

The main aim of this paper is to study the generic properties of finitely
presented monoids and semigroups, and hence to understand the generic-case
complexity of uniform decision problems for monoids and semigroups.
Our main results show that, with respect to a number of very natural
stratifications, the generic\footnote{For brevity, we use statements
such as ``the generic
$X$ has property $Y$'' as shorthand for ``there is a generic subset of
the set of $X$'s, every member of which has property $Y$''. Of course
the generic $X$ truly ``exists'' only in the case that a single
isomorphism type forms a generic subset of $X$'s; in this case the
isomorphism type has all the ascribed properties, so the terminology
is unambiguous!} finite monoid presentation (over a given alphabet and
with a given number of generators) satisfies \textit{small overlap
conditions} in the sense introduced by
Remmers \cite{Remmers71,Remmers80} (see also \cite{Higgins92}).
Small overlap conditions are natural semigroup-theoretic analogues of the
\textit{small overlap conditions} extensively used by combinatorial group
theorists, and so our main result can be viewed as loosely analogous
(although our objectives and hence our formalism are rather different)
to the well-known fact, first asserted by Gromov \cite{Gromov87} and proved
in detail by Ol'shanskii \cite{Olshanskii92}, that the generic finitely presented group is word
hyperbolic.

These results immediately tell us a great deal about the algebraic structure
of the generic monoid. For example, we
learn that
it is \textit{$\mathscr{J}$-trivial}, and hence torsion-free with no
non-trivial subgroups. Even more important, by recent results of the
author \cite{K_smallover}, the uniform word problem for such presentations
is solvable in (worst-case) time linear in the words lengths and quadratic
in the presentation size. Since it can be checked in (worst-case) quadratic
time whether a presentation satisfies a small overlap condition, it follows that the uniform
word problem for finitely presented monoids is generically solvable in
(worst-case) time linear in the word lengths and quadratic in the
presentation size. All of these results apply equally to semigroups without
identity elements.

As already remarked, generic-case complexity has been developed by group theorists,
and the literature is largely concerned with applications to advanced group
theory; as a result, much of it is not readily accessible to non-algebraists.
An additional objective of this article is to
provide a gentle (although by no means comprehensive) introduction to generic
sets and properties, and generic-case complexity,
in a form fully intelligible to the reader without a specialist algebraic
background. Monoid presentations are combinatorially simpler objects than
group presentations, and most of our proofs are of an elementary combinatorial
nature which should allow them double as detailed worked examples to give the
reader a feel for the theory of generic-case complexity. The few places where we
resort to more advanced algebraic notions are clearly delimited and
self-contained, so that the bulk of the article can be understood without
following these parts in detail.

In addition to this introduction, this article comprises four sections.
Section~\ref{sec_generic} provides a gentle introduction to generic sets
and generic-case complexity. In Section~\ref{sec_monoids} we prove our main
results about generic monoids and semigroups with respect to certain
stratifications.
In Section~\ref{sec_equivalence} we prove some technical results regarding
the relationships between different stratifications; these may be of some
independent interest; these are applied to show that our results about
generic monoids apply regardless of which of several natural stratifications
is chosen. Finally, Section~\ref{sec_consequences} explores the
consequences of our characterisations of generic monoids and semigroups,
including the fact that the uniform word problems for finitely presented
monoids and semigroups are generically solvable in time quadratic in the
presentation lengths and linear in the word lengths.

\section{Generic Properties and Generic-case Complexity}\label{sec_generic}

In this section we provide a brief introduction to generic sets and generic
complexity. A more comprehensive treatment can be found in \cite{Gilman07}.
Our aim is to make the paper accessible to as wide an audience as possible,
and so we endeavour to keep mathematical prerequisites to a minimum. However,
we cannot avoid assuming some elementary familiarity with the theory of sets
and sequences.

Let $S$ be a countably infinite set. A \textit{stratification} of
$S$ is an infinite sequence $S_1, S_2, \dots, S_n, \dots$ of finite subsets
of $S$ whose union is $S$. The computationally-orientated reader may like to
bear in mind the example where $S$ is the instance space for some problem,
and $S_n$ is the set of instances of size $n$ for some suitable notion of
size; however, we caution that in general the subsets $S_n$ need not be
disjoint. We call the stratification \textit{spherical} if the sets $S_n$ are
pairwise disjoint ($S_i \cap S_j = \emptyset$ for all $i \neq j$) and at
the other extreme \textit{ascending} if they form an ascending sequence
under containment ($S_i \subseteq S_j$ for all $i < j$).

Now let $X$ be a subset of $S$. We say that $X$ is \textit{generic} (with respect to the
given stratification) if
\begin{equation}\label{eq_generic}
\lim_{n \to \infty} \frac{|X \cap S_n|}{|S_n|} = 1
\end{equation}
The subset $X$ is called \textit{negligible} if $S \setminus X$ is
generic, or equivalently, if the limit on the left-hand-side of
\eqref{eq_generic} is defined and equal to $0$.
Intuitively, $X$ is generic if the probability that an instance of size
$n$, chosen uniformly at random, lies in $X$ can be made arbitrarily close
to $1$ by choosing large enough $n$.

Note that, for any given set $X$, the limit on the left-hand-side of
\eqref{eq_generic} may not be defined, and indeed for almost any
stratification it is easy to construct a set $X$ for which it is not.
The function
$$X \mapsto \lim_{n \to \infty} \frac{|X \cap S_n|}{|S_n|}$$
is a finitely additive probability measure defined on those subsets of $X$ for
which the limit converges, but it is typically \textit{not} a measure in the
usual sense, since it lacks countable additivity. This fact is no accident:
a countably additive measure on a countable set cannot assign $0$-measure to
all the singletons, but as we noted in the introduction, a key feature of the
generic approach is that single instances are regarded as negligible.
Nonetheless, the intuition that the generic sets are those of ``full measure''
can be helpful, and they satisfy many of the elementary properties of such
sets. In particular, the reader can easily verify that if $X$ is generic
and $X \subseteq Y$ then $Y$ is generic, while if $X$ and $Y$ are both generic
then $X \cap Y$ is generic. Obvious dual statements hold for negligible sets.

Notice that, in our initial definition of generic sets, we have placed no
requirements on the \textit{rate} of convergence of the left-hand-side of
\eqref{eq_generic}. Genericity is an asymptotic property, and if convergence is very slow then
the asymptotic behaviour may not be reflected in ``practical sized'' instances.
We call a set $X$ \textit{superpolynomially} generic/negligible if the
appropriate limit converges faster than $1/n^p$ for \textit{every} $p \in \mathbb{N}$,
and \textit{exponentially} generic/negligible if it converges faster than $p^n$ for
for {some} $p \in (0, 1]$.  (In the literature some authors use the term
``strongly generic'' for what we have called exponentially generic sets, while
some use ``strongly generic'' to mean superpolynomially generic and
``supergeneric'' to mean exponentially generic. To avoid confusion, we shall
avoid these terms in favour of less concise but more descriptive ones.)

We now turn our attention to the application of generic sets in
computational complexity. This requires us to consider explicitly
not just abstract algorithmic problems, but also also stratifications of
instance spaces. We define a \textit{stratified problem} to be an
algorithmic decision problem equipped with a stratification on its instance
space. (We shall restrict our attention here to decision problems,
but analogous definitions can be made for more general computational problems.)

Of course traditional complexity theory is implicitly concerned with
stratified problems: to study the asymptotic complexity of a problem one
requires a notion of the \textit{size} of each member of the instance space
$S$. As we have already remarked, this automatically induces a stratification
given by setting $S_n$ to be the set of all instances of size $n$. We call
this the \textit{input size stratification} for the problem. However, the
dependence on stratification is much tighter in generic
complexity theory than it is in traditional complexity theory -- many
authors discussing traditional complexity of algorithmic problems prefer
to avoid detailed discussion of data encoding and hence of exact instance sizes; this is
entirely reasonable since traditional complexity classes are largely
insensitive to minor encoding issues. But for generic-case complexity, these
issues can make a very big difference.

Note also that, while the input size stratification is a natural,
canonical one to associate to any algorithmic problem, it is only one
of many possible stratifications, and may not be the
appropriate one for any given application. The ideal is rather to find a
stratification which reflects the empirical distribution of problem
instances, that is, the frequency with which they arise in practice in
a particular application, and there is often no reason to suppose
that this is strongly correlated with size.

Now let $\mathbf{C}$ be any class of decision problems (typically a
complexity class of some kind). We say that a stratified problem
$\mathscr{P}$ is \textit{generically in $\mathbf{C}$} if there exists
a generic subset $Y$ of the instance space such that
\begin{itemize}
\item[(i)] the membership problem for $Y$ lies in $\mathbf{C}$; and
\item[(ii)] the problem $\mathscr{P}$ restricted to $Y$ lies in $\mathbf{C}$.
\end{itemize}
Intuitively, a stratified decision problem is generically in $\mathbf{C}$
if the decision problem admits a partial algorithm (that is, an algorithm
which outputs ``yes'', ``no'' or ``don't know'', and which in the former two
cases is always correct) in $\mathbf{C}$, such that the probability of a
``don't know'' is negligible. We write $\mathbf{GenC}$ for the class of
all stratified problems generically in $\mathbf{C}$.

Obvious examples are the class $\mathbf{GenP}$ of \textit{generically polynomial-time stratified
problems} and \textbf{GenNP} of \textit{generically non-deterministic
polynomial-time stratified problems}.
Another interesting example is the class $\mathbf{GenBPP}$, which consists of stratified
problems admitting a randomised polynomial-time algorithm with probabability
of error uniformly bounded away from $1/2$ for every instance in some generic
subset whose membership problem also lies in $\mathbf{BPP}$.

\section{Generic Monoid Presentations}\label{sec_monoids}

In this section we study the generic properties of finite monoid
presentations. We begin with some basic definitions.

Let $A$ be a finite alphabet (set of symbols). A \textit{word} over $A$
is a finite sequence of zero or more elements from $A$. The set of all words
over $A$ is denoted $A^*$; under the operation of \textit{concatenation}
it forms a monoid, called the \textit{free monoid} on $A$. The length of
a word $w \in A^*$ is denoted $|w|$. The unique \textit{empty word} of length
$0$ is denoted $\epsilon$; it forms the identity element of the monoid
$A^*$. The set $A^+ \setminus \lbrace \epsilon \rbrace$ of non-empty
words forms a subsemigroup of $A^*$, called the \textit{free semigroup
on $A$}.

A finite monoid presentation $\langle A \mid R \rangle$ consists of a finite
alphabet $A$, together with a finite sequence $R \subseteq A^* \times A^*$ of
ordered pairs of words\footnote{The reader may think it more natural to
consider a \textit{set} of \textit{unordered} pairs, but the definition we
use simplifies the combinatorics in our analysis, and Theorem~\ref{thm_ordervsunorder}
will show that it makes no difference to the end results.}.
 We say that $u, v \in A^*$ are
\textit{one-step equivalent} if $u = axb$ and $v = ayb$ for some
possibly empty words $a, b \in A^*$ and relation $(x,y) \in R$ or $(y,x) \in R$.
We say that $u$ and $v$ are \textit{equivalent}, and write $u \equiv_R v$
or just $u \equiv v$, if there is a finite sequence of words beginning
with $u$ and ending with $v$, each term of which but the last is one-step
equivalent to its successor. Equivalence is clearly an equivalence relation;
in fact it is the least equivalence relation containing $R$ and compatible
with the multiplication on $R$. The equivalence classes form a monoid
with multiplication well-defined by $[u]_\equiv [v]_\equiv = [uv]_\equiv$;
this is called the \textit{monoid presented} by the presentation.

The \textit{word problem} for a (fixed) monoid presentation
$\langle A \mid R \rangle$ is the algorithmic problem of, given as
input two words $u, v \in A^*$, deciding whether $u \equiv_R v$. The
\textit{uniform word problem for finitely presented monoids} is the
algorithmic problem of, given as input a monoid presentation
$\langle A \mid R \rangle$ and two words $u, v \in A^*$, deciding
whether $u \equiv_R v$. It is well-known that there exist finite
monoid presentations which the word problem is undecidable, and hence
that the uniform word problem for finitely presented monoids is
undecidable \cite{Markov47,Post47}. More generally, if $\mathscr{C}$ is a class
of finite monoid presentations, then the uniform word problem for
$\mathscr{C}$ monoids is the algorithmic problem of, given as input a
monoid presentation $\langle A \mid R \rangle$ in $\mathscr{C}$ and
two words $u, v \in A^*$, deciding whether $u \equiv_R v$.

Now suppose we have a fixed monoid presentation $\langle A \mid R \rangle$.
A \textit{relation word} is a word which appears as one side of a relation
in $R$. A \textit{piece} is a word which appears more than once as a factor
in the relations, either as a factor of two different relation words, or as
a factor of the same relation word in two different (but possibly overlapping)
places. Let $m \in \mathbb{N}$ be a positive integer. The presentation is said to 
\textit{satisfy $C(m)$} if no relation word can be written as a product of
\textit{strictly fewer than $m$} pieces. Thus $C(1)$ says that no relation
word is empty; $C(2)$ says that no relation word is a factor of another.

Definitions corresponding to all of those above can also be made for
semigroups (without necessarily an identity element), by taking $A^+$ in
place of $A^*$ (in all places except the definition of one-step equivalence,
where $a$ and $b$ must still be allowed to be empty).

Now fix an alphabet $A$. To study generic properties of $k$-relation presentations
over $A$, we need a stratification on the (countable) set of all such. There
are two obvious ways to define the size of a presentation, and hence two natural stratifications of
the $A$-generated $k$-relation presentations.
Firstly, one can
take the size of the presentation to be the \textit{sum
length of the relation words}; this gives rise to the \textit{sum length
stratification} of presentations. Alternatively, one can define the
size to be the \textit{length of the longest relation word}; this results in
the \textit{maximum relation stratification}. Which choice is most
natural depends on the application. For example, the sum length of
a presentation is a good approximation to the space required to encode the
presentation in the obvious way, and hence for computational applications
seems the most natural. Intuitively, the sum length
stratification lends greater weight to uneven distributions of the relation
word lengths within a presentation; in particular, it results in a greater
frequency of short words, which makes it seem less likely that small overlap
conditions will hold. Nevertheless, it transpires that our main results
hold for both stratifications, which may be regarded as some evidence of
their ``robustness''.

We emphasise that we are attempting here to stratify only the set of
$A$-generated, $k$-relation semigroup presentations, where the alphabet $A$
and set of relations $k$ are fixed. There are, of course, also natural
stratifications across \textit{all} $A$-generated semigroup presentations,
allowing the number of relations to vary. These typically lead to a
high frequency of ``short'' relation words, which means that
small overlap type conditions do not hold generically. However, it seems
likely that, for at least some natural stratifications of this type, the
word problem remains generically solvable for other reasons. This
interesting issue will be studied further in a subsequent paper.

We shall need a couple of elementary definitions from combinatorics. Let $n$
and $k$ be non-negative integers. Recall that a \textit{composition of $n$
into $k$} is an ordered $k$-tuple of positive integers which sum to $n$,
while a \textit{weak composition of $n$ into $k$} is an ordered $k$-tuple
of non-negative integers which sum to $n$.

Having fixed the alphabet $A$, a $k$-relation monoid presentation of sum
length $n$ is uniquely determined by its sequence of relation words; this
in turn is uniquely determined by the concatenation in order of those words
(a word in $A^n$) and the lengths of those words (a weak composition of
$n$ into $2k$, called the \textit{shape} of the presentation). Thus, $k$-relation monoid presentations of sum relation
length $n$ are in a bijective correspondence with ordered pairs whose
first component is a word of length $n$, and whose second component is a
weak composition of $n$ into $2k$.

We shall need the following simple combinatorial lemma.
\begin{lemma}\label{lemma_overlap}
Let $A$ be a finite alphabet and $c$ and $p$ be positive integers. The
number of distinct words of length $c$ which admit factorisations as
$x_1vy_1$ and as $x_2vy_2$ for some $x_1, x_2, y_1, y_2, v \in A^*$ with
$|v| \geq p$ and $x_1 \neq x_2$ is bounded above by $c^2 |A|^{c-p}$.
\end{lemma}
\begin{proof}
Clearly if a word admits such factorisations, then it admits such
factorisations with $|v| = p$, so we may count only those words which
admit such factorisations with $|v| = p$.

We claim, having fixed $A$, $c$ and $p$, any such word is uniquely determined by $x_1$, $y_1$ and the
length of $x_2$. Clearly, there are fewer than $c^2$ ways to choose the
lengths of $x_1$ and $x_2$; doing so also fixes the length of $y_1$, since
we must have
$$|x_1| + |v| + |y_1| = |x_1| + p + |y_1| = c.$$
Now there are at most
$$|A|^{|x_1| + |y_1|} \ = \ |A|^{c - |v|} \ = \ |A|^{c-p}$$
ways to choose the words $x_1$ and $y_1$ with the given lengths, so
proving the claim will suffice to prove the lemma.

Since $x_1$ and $x_2$ are distinct prefixes of the same word, their lengths
cannot be equal.
Suppose first that $x_1$ is longer than $x_2$ and write
$v = v^{(1)} \dots v^{(|v|)}$ and $x_1 = x_1^{(1)} \dots x_1^{(|x_1|)}$ with each
$v^{(i)}$ and $x_1^{(i)}$ in $A$. Then since $x_1 v y_1 = x_2 v y_2$ we have
$$v^{(i)} = \begin{cases}
      x_1^{(|x_2| + i)}        &\text{ for } 1 \leq i \leq |x_1| - |x_2| \\
      v^{(i - |x_1| + |x_2|)}  &\text{ for } |x_1| - |x_2| < i \leq |v|
  \end{cases}
$$
from which the claim follows.

If, on the other hand, $x_1$ is shorter than $x_2$ then we use the lengths
of $v$ and $x_2$ to deduce the length of $y_2$, whereupon a
symmetric argument suffices to complete the proof.
\end{proof}

\begin{proposition}\label{prop_piecelength_fixedshape}
Let $A$ be a finite alphabet, and $n$ and $r$ be positive integers, and
fix a weak composition $\sigma$ of $n$ (into any number).
Then the proportion of presentations of shape $\sigma$ which have a
piece of length $r$ or more is bounded above by $n^2 |A|^{-r}$.
\end{proposition}
\begin{proof}
The set of presentations over $A$ of shape $\sigma$ is in 1:1 correspondance
with the set $A^n$ via the map which takes each presentation to the
concatenation, in the obvious order, of its relation words. If the
presentation has a piece of length $r$ or more then the corresponding
word will feature that piece as a factor in at least two different places.
By Lemma~\ref{lemma_overlap} it follows that the number of presentations
with a piece of length $r$ or more is bounded above by
$n^2 |A|^{n-r}$. The total number of such presentations in $|A|^n$, so
the proportion of presentations with the desired property is bounded above
by $n^2 |A|^{-r}$
as required.
\end{proof}

\begin{corollary}\label{cor_smalloverlap_fixedshape}
Let $A$ be a finite alphabet and $k$, $n$, $m$ and $K$ be positive integers
with $m \geq 2$,
and fix an weak composition
$\sigma$ of $n$ into $2k$ such that no block has size less than. Then
the proportion of presentations with alphabet $A$ and shape $\sigma$ which
do not satisfy $C(m)$ is bounded above by
$$\frac{n^2}{|A|^{K/(m-1)}}.$$
\end{corollary}
\begin{proof}
If a presentation fails to satisfy $C(m)$ then some relation word can be
written as a product of $m-1$ pieces. By assumption this relation word must
have length at least $K$, so one of the pieces must have length at least
$K / (m-1)$. The result is now immediately from
Proposition~\ref{prop_piecelength_fixedshape}.
\end{proof}

Before proving the first of our main theorems, we will need an elementary
combinatorial result concerning weak compositions; this will serve to bound
the proportion of presentations which feature a ``short'' relation word.

\begin{lemma}\label{lemma_fewshortwords_sumlength}
Let $k$ be an integer, and $f : \N \to \N$ be a function such that
$f(n) / n$ tends to zero as $n$ tends to infinity. Then the proportion of weak
compositions of $n$ into $k$ which feature a block of size $f(n)$ or less
tends to zero as $n$ tends to infinity.
\end{lemma}
\begin{proof}
It is well-known and easy to prove (see, for example, \cite[Theorem~5.2]{Bona02}) that the number of weak compositions
of $n$ into $k$ is given by
$$C'_k(n) \ = \ \frac{(n+k-1)!}{n! \ (k-1)!}$$
Clearly, every partition of $n$ into $k$ featuring a block of size
$f(n)$ or less can be obtained by refining a partition of $n$ into
$k-1$, with the extra decomposition in one of $k (f(n)+1)$ places.
Thus, the number of such partitions is bounded above by
$$k \ (f(n) + 1) \ C'_{k-1}(n) \ = \ 
k \ (f(n) + 1) \ \frac{(n+k-2)!}{n! \ (k-2)!}$$
Hence, the proportion of such partitions amongst all weak compositions
of $n$ into $k$ is bounded above by
\begin{align*}
\frac{k \ (f(n) + 1) \ C'_{k-1}(n)}{C'_k(n)} \ 
&= \ \frac{k \ (f(n) + 1) \ (n+k-2)! \ n! \ (k-1) !}{ (n+k-1)! \ n! \ (k-2)!} \\
&= \ \frac{k (k-1) \ (f(n) + 1)}{n+k-1} \\
&= \ k (k-1) \left( \frac{f(n)}{n+k-1} + \frac{1}{n+k-1} \right) \\
&\leq \ k (k-1) \left( \frac{f(n)}{n} + \frac{1}{n+k-1} \right)
\end{align*}
which clearly tends to zero as $n$ tends to infinity.
\end{proof}

We are now ready to prove our main theorem for the sum relation
length stratification.
\begin{theorem}\label{thm_generic_sumlength}
Let $A$ be an alphabet of size at least $2$, and $k$ and $m$ be positive
integers. Then the set of $A$-generated, $k$-relation monoid presentations
which satisfy the condition $C(m)$ is generic with respect to the sum
length stratification.
\end{theorem}
\begin{proof}
Since $C(2)$ implies $C(1)$, we may clearly assume without loss of
generality that $m \geq 2$. We need to show that the proportion of $A$-generated $k$-relation monoid
presentations of length $n$ which fail to satisfy $C(m)$ tends to zero
as $n$ tends to infinity.

For each $n$, let $P_n$ be the set of all weak compositions of $n$ into $k$,
let $Q_n$ be the set of weak compositions of $n$ into $k$ featuring a
block of size $3(m-1) \log_{|A|} n$ or less, and let $R_n = P_n \setminus Q_n$.
By an application of Lemma~\ref{lemma_fewshortwords_sumlength}, with the function
$f : \N \to \N$ given by $f(n) = 3(m-1) \log_{|A|} n$, we see that the proportion
$|Q_n| / |P_n|$ tends to $0$ as $n$ tends to infinity.

For each weak composition $\sigma$, let $x_\sigma$ be the proportion of
presentations of shape $\sigma$ which fail to satisfy $C(m)$. Note that
by Corollary~\ref{cor_smalloverlap_fixedshape} we have
$$x_\sigma \ \leq \ \frac{n^2}{|A|^{K_\sigma / (m-1)}}$$
where $K_\sigma$ denotes the smallest block size in $\sigma$.
For each
fixed $n$, there are clearly equally many ($|A|^n$ to be precise)
presentations of each shape, so the
proportion of presentations of length $n$ failing to satisfy $C(m)$ is just the
average over shapes $\sigma$ of $x_\sigma$, that is:
\begin{align*}
\frac{1}{|P_n|} \left( \sum_{\sigma \in P_n} x_\sigma \right)
\ &= \ \frac{1}{|P_n|} \left( \sum_{\sigma \in Q_n} x_\sigma \right) + \frac{1}{|P_n|}  \left( \sum_{\sigma \in R_n} x_\sigma \right) \\
&\leq \ \frac{1}{|P_n|} \left( \sum_{\sigma \in Q_n} 1 \right) + \frac{1}{|P_n|}  \left( \sum_{\sigma \in R_n} \frac{n^2}{|A|^{K_\sigma/(m-1)}} \right) \\
&= \ \frac{|Q_n|}{|P_n|} + \frac{1}{|P_n|} \left( \sum_{\sigma \in R_n} \frac{n^2}{|A|^{K_\sigma / (m-1)}} \right).
\end{align*}
We have already observed that $|Q_n| / |P_n|$ tends to zero as $n$
tends to infinity. Moreover, by the definition of $R_n$ we have
$K_\sigma > 3(m-1) \log_{|A|} n$ for all $\sigma \in R_n$ so that
\begin{align*}
\frac{1}{|P_n|} \sum_{\sigma \in R_n} \frac{n^2}{|A|^{K_\sigma / (m-1)}} \ 
&\leq \ \frac{1}{|P_n|} \sum_{\sigma \in R_n} \frac{n^2}{|A|^{(3(m-1) \log_{|A|} n) / (m-1)}} \\
&= \ \frac{|R_n|}{|P_n|} \frac{n^2}{|A|^{(3(m-1) \log_{|A|} n) / (m-1)}} \\
&= \ \frac{|R_n|}{|P_n|} \frac{n^2}{|A|^{\log_{|A|} (n^3)}} \\
&\leq \ \frac{n^2}{n^3}
\end{align*}
which tends to zero as required.
\end{proof}

An analysis of the proof shows, approximately speaking, that the proportion
of presentations of $A$ failing to satisfy any given small overlap
condition goes to zero like $(\log_{|A|} n) / n$, which for practical
purposes may be rather slow. The barrier to showing a faster convergence is the proportion
of presentations featuring a ``short'' relation word ($|Q_n|/|P_n|$ in the
notation of the proof); this proportion really does seem to decrease very
slowly, suggesting that for the sum length stratification, fast
convergence to small overlap conditions is not possible. To obtain statements
about the ``superpolynomially generic monoid'' or ``exponentially generic
monoid'' with respect to the sum length stratification, one would require
arguments which take detailed account of the ``short'' relation words.

Our next task is to prove that an equivalent result holds for the maximum
length stratification. We begin with an analogue of
Lemma~\ref{lemma_fewshortwords_sumlength}, which will show that the frequency
of presentations featuring a ``small'' relation word is again negligible.
This time, because the number of presentations of each shape of maxmimum
length $k$ is not fixed, we must reason directly with presentations rather
than just shapes. Having taken account of this, the result is easier and,
as one might expect given our remarks above on the relative frequency of
``short'' relation words in this stratification, stronger.

\begin{lemma}\label{lemma_fewshortwords_maxlength}
Let $A$ be an alphabet of size at least $2$, $k$ be a non-negative integer, and
$f : \N \to \N$ be a function such that $n - f(n)$ tends to infinity as $n$ tends to infinity. Then the proportion of
$A$-generated $k$-relation presentations of maximum relation word length $n$ which feature a relation
word of length $f(n)$ or less tends to zero as $n$ tends to infinity. Moreover,
if there exists a constant $p > 0$ such that $n - f(n) > p n$ for sufficiently
large $n$ then the given proportion tends to zero exponentially fast.
\end{lemma}
\begin{proof}
Let $X_n$ be the set of all presentations over $A$ of maximum relation
length $n$, let $Y_n$ be the presentations in $X_n$ which have a relation
word of
length $f(n)$ or less, and let $Z_n = X_n \setminus Y_n$. The quantity we seek is thus the limit
as $n$ tends to infinity of $|Y_n| / |X_n|$. Let $I = \lbrace 1, \dots, 2 k \rbrace$ and define a
map $\sigma$ from $I \times X_n$ to the set of all presentations $k$-relation
presentations of $A$, 
which takes $(i, P)$ to the presentation obtained from $P$ by removing
$n-f(n)$ characters from the end of the $i$th relation word, or replacing
this relation word with the empty word if its length is less than $n-f(n)$.

We claim that under the map $\sigma$, every presentation in $Y_n$ has
at least $|A|^{n-f(n)}$ pre-images in $I \times X_n$.
Indeed, if $Q \in Y_n$ then $Q$ has some relation word (say the $j$th)
of length less than $f(n)$, say length $p$. Now for each of $|A|^{n-f(n)}$ words $w \in A^{n-f(n)}$ we can obtain from $Q$ a
presentation $P_w \in X_n$ by appending $w$ to the end of the $j$th
relation word, and it is easily seen $\sigma(j, P_w) = Q$ for all such $w$.

Thus, we have $2k |X_n| = | I \times X_n | \geq |A|^{n - f(n)} |Y_n|$, and so
\begin{align*}
\frac{|Y_n|}{|X_n|} \ &\leq \ \frac{2k}{|A|^{n-f(n)}}.
\end{align*}
Since $n-f(n)$ tends to infinity with $n$, this clearly tends to zero.
If moreover $p > 0$ is such that $n - f(n) \geq p n$ for $n$ sufficiently
large then we have
$$\frac{|Y_n|}{|X_n|} \ \leq \ \frac{2k}{|A|^{pn}}$$
so that the given quantity tends to zero exponentially fast.
\end{proof}

\begin{corollary}\label{cor_fewshortwords_maxlength}
Let $A$ be an alphabet of size at least $2$, $k$ be a non-negative integer,
and $c$ a constant with $0 < c < 1$. Then the proportion of
$A$-generated, $k$-relation presentations of maximum relation word length $n$ which feature a relation
word of length $cn$ tends to zero exponentially fast as $n$ tends to infinity.
\end{corollary}
\begin{proof}
Define $f : \mathbb{N} \to \mathbb{N}$ by $f(n) = cn$, and choose $p$ with
$$0 < p < 1-c.$$
Then $n - f(n) = (1-c)n > pn$ for all $n$, so the result
follows from Lemma~\ref{lemma_fewshortwords_maxlength}.
\end{proof}

We are now ready to prove our main result for the maximum length
stratification.
\begin{theorem}\label{thm_generic_maxlength}
Let $A$ be an alphabet of size at least $2$, and let $k$ and $m$ be positive integers.
Then the set of $A$-generated, $k$-relation monoid presentations which
satisfy $C(m)$ is exponentially generic with respect to the maximum
length stratification.
\end{theorem}
\begin{proof}
The structure of the proof is essentially the same as that for
Theorem~\ref{thm_generic_maxlength}, but it is slightly complicated
by the fact that the number of presentations of each shape for a given
maximum relation word $n$ is not fixed. In addition, we must to show
that the rate of convergence is exponential. Once again, we assume
without loss of generality that $m \geq 2$.

Let $C_n$ be the total number of presentations over $A$ of maximum relation
word length $n$.
Let $P_n$ be the set of all weak compositions of any integer into $2k$
with largest block size $n$.
Choose $d$ with $0 < d < 1$ and let $Q_n$ be the set of all shapes in $P_n$
with a word of length $dn$ or less.
Let $R_n = P_n \setminus Q_n$.
For each weak composition $\sigma \in P_n$, let $c_\sigma$ be the total
number of presentations of shape $\sigma$, and let $x_\sigma$ be the
proportion of presentations of shape $\sigma$ which fail to satisfy
$C(m)$. For each shape $\sigma$, by Corollary~\ref{cor_smalloverlap_fixedshape} we
have
$$x_\sigma \ \leq \ \frac{(n_\sigma)^2}{|A|^{K_\sigma / (m-1)}}$$
where $n_\sigma$ is the total size of $\sigma$ (that is, the sum of the
block sizes of $\sigma$, or the \textit{sum} relation word length of a presentation
of shape $\sigma$), and $K_\sigma$
is the smallest block size in $\sigma$. But $\sigma$ has $2k$ blocks, none
of which is larger than $n$, so we must have
$n_\sigma \leq 2kn$, so that
$$x_\sigma \ \leq \ \frac{(2kn)^2}{|A|^{K_\sigma/(m-1)}} \ = \ \frac{4 \ k^2 \ n^2}{|A|^{K_\sigma/(m-1)}}.$$
Now the proportion we seek is given by
\begin{align*}
\frac{1}{C_n} \left( \sum_{\sigma \in P_n} c_\sigma x_\sigma \right)
\ &= \ \frac{1}{C_n} \left( \sum_{\sigma \in Q_n} c_\sigma x_\sigma \right) + \frac{1}{C_n}  \left( \sum_{\sigma \in R_n} c_\sigma x_\sigma \right) \\
&\leq \ \frac{1}{C_n} \left( \sum_{\sigma \in Q_n} c_\sigma \right) +
 \frac{1}{C_n}  \left( \sum_{\sigma \in R_n} c_\sigma \frac{4 k^2 n^2}{|A|^{K_\sigma/(m-1)}} \right).
\end{align*}
The first term in the last line is the proportion of presentations
featuring a relation word of length $dn$ or less; by Corollary~\ref{cor_fewshortwords_maxlength},
this tends to zero exponentially fast. Considering now the second term, by
the definition of $R_n$ we have that $K_\sigma > dn$ for all $\sigma \in R_n$
so that
\begin{align*}
\frac{1}{C_n} \sum_{\sigma \in R_n} c_\sigma \frac{4 k^2 n^2}{|A|^{K_\sigma / (m-1)}} \ 
&\leq \ \frac{1}{C_n} \sum_{\sigma \in R_n} c_\sigma \frac{4 k^2 n^2}{|A|^{dn / (m-1)}} \\
&= \ \left( \frac{4 k^2 n^2}{|A|^{dn/(m-1)}} \right) \ \left( \frac{\sum_{\sigma \in R_n} c_\sigma}{C_n} \right) \\
&\leq \ \frac{4 k^2 n^2}{(|A|^{d/(m-1)})^n}.
\end{align*}
which since $|A| \geq 2$ and $d > 0$ clearly tends to zero exponentially fast.
\end{proof}

\section{Equivalence of Stratifications}\label{sec_equivalence}

It often happens that two stratifications (on the same set, or on related
sets) are closely related, so that knowledge of the generic sets with respect
to one yields corresponding information about the generic sets with respect
to the other. In this section we establish some technical conditions under
which this holds, and use this to extend many of our earlier results to
additional natural stratifications.

First, we consider the relationship between spherical and ascending
stratifications. So far, we have seen examples only of spherical
stratifications of instance spaces, but to each such stratification is
associated an equally natural ascending stratification, the sets in the
latter being unions of the sets in the former. The following proposition,
which was first observed in \cite{Gilman07}
to be an easy consequence of the Stolz-Cesaro Theorem, says
that the generic sets are independent of which of these stratifications
is used (see \cite{Gilman07}
for a more detailed explanation).
\begin{proposition}{\cite[Lemma~3.2]{Gilman07}}
Let $S_n$ be a spherical stratification of a set $S$. Define a new
stratification on $S$ by
$$B_n \ = \ \bigcup_{j = 1}^{n} S_j.$$
Then any set $X \subseteq S$ is generic with respect to the stratification
$S_n$ if and only if it is generic with respect to the stratification $B_n$.
\end{proposition}

We shall need the following elementary proposition, which essentially
says that the restriction of a stratification to a generic set preserves
generic sets.
\begin{lemma}\label{lemma_genericsubsetofgeneric}
Let $X$ be a stratified set, and $X'$ a generic subset of $X$. Then for
any $P \subseteq X$ we have
$$\lim_{n \to \infty} \frac{|P \cap X_n|}{|X_n|} = \lim_{n \to \infty} \frac{|P \cap X_n \cap X'|}{|X_n \cap X'|}.$$
\end{lemma}
\begin{proof}
First notice that, since $X'$ is generic, we have
\begin{equation}\label{eq_1}
\lim_{n \to \infty} \frac{|P \cap X_n \cap (X \setminus X')|}{|X_n|}
\ = \ \lim_{n \to \infty} \frac{|(X \setminus X') \cap X_n|}{|X_n|} \ = \ 0 
\end{equation}
Now
\begin{align*}
\lim_{n \to \infty} \frac{|P \cap X_n \cap X'|}{|X_n \cap X'|} \ 
&= \ \lim_{n \to \infty} \frac{|P \cap X_n \cap X'|}{|X_n|} \ \frac{|X_n|}{|X_n \cap X'|} \\
&= \ \left( \lim_{n \to \infty} \frac{|P \cap X_n \cap X'|}{|X_n|} \right) \ \left( \lim_{n \to \infty} \frac{|X_n \cap X'|}{|X_n|} \right)^{-1} \\
&= \ \left( \lim_{n \to \infty} \frac{|P \cap X_n \cap X'|}{|X_n|} \right) \ 1^{-1} \ \ \  \text{(since $X'$ is generic)} \\
&= \ \left( \lim_{n \to \infty} \frac{|P \cap X_n \cap X'|}{|X_n|} \right) \ + \ 0 \\
&= \ \left( \lim_{n \to \infty} \frac{|P \cap X_n \cap X'|}{|X_n|} \right) \ + \left( \ \lim_{n \to \infty} \frac{|P \cap X_n \cap (X \setminus X')}{|X_n|} \right) \text{(by \eqref{eq_1})} \\
&= \ \lim_{n \to \infty} \frac{|P \cap X_n \cap X'|}{|X_n|} + \frac{|P \cap X_n \cap (X \setminus X')}{|X_n|} \\
&= \ \lim_{n \to \infty} \frac{|P \cap X_n|}{|X_n|}
\end{align*}
as required.
\end{proof}

Next, we introduce a very useful sufficient condition for a map between
stratified sets to preserve generic sets. To do so, we need some terminology.
Let $X$ and $Y$ be stratified sets, $X' \subseteq X$ and $Y' \subseteq Y$,
and $f : X' \to Y'$ a map. Then $f$ is called \textit{stratification-preserving} if for
every $x \in X'$ and $n \in \mathbb{N}$ we have $x \in X_n$ if and only if
$f(x) \in Y_n$. If $P \subseteq X$ then $f$ is said to \textit{respect $P$}
if $f(P \cap X')$ and $f((X \setminus P) \cap X')$ are disjoint, that is, 
if whenever $x_1, x_2 \in X'$ are such that $f(x_1) = f(x_2)$ we have
either $x_1, x_2 \in P$ or $x_1, x_2 \notin P$.
 Recall that the \textit{fibre size} of $f$ at a point
$y \in Y'$ is the cardinality of the set of elements $x \in X'$ such that
$f(x) = y$. The map $f$ is called \textit{bounded-to-one} if there is a
finite upper bound on its fibre sizes.

\begin{theorem}\label{thm_equivalence}
Let $X$ and $Y$ be stratified sets, $X' \subseteq X$ and $Y' \subseteq Y$
be generic subsets of $X$ and $Y$ respectively, $d \in \mathbb{N}$ and
$f : X' \to Y'$ a surjective, stratification-preserving map, such that for
every $n \in \mathbb{N}$ there exists $k_n \in \mathbb{N}$ such that the
fibre sizes of $f$ at points in $X_n \cap X'$ all lie between $k_n$ and
$d k_n$. Then for any set $P \subseteq X$ we have
\begin{itemize}
\item[(i)]
$$\frac{1}{d} \lim_{n\to\infty} \frac{|f(P \cap X') \cap Y_n|}{|Y_n|} \ \leq \ 
\lim_{n \to \infty} \frac{|P \cap X_n|}{|X_n|} \ \leq \ 
d \lim_{n\to\infty} \frac{|f(P \cap X') \cap Y_n|}{|Y_n|}$$
wherever both limits are defined;

\item[(ii)]
$$\frac{1}{d} \lim_{n\to\infty} \frac{|P \cap X_n|}{|X_n|} \ \leq \ 
\lim_{n \to \infty} \frac{|f(P \cap X') \cap Y_n|}{|Y_n|} \ \leq \ 
d \lim_{n\to\infty} \frac{|P \cap X_n|}{|X_n|}$$
wherever both limits are defined;

\item[(iii)] $P$ is negligible in $X$ if and only if
$f(P \cap X')$ is negligible in $Y$;

\item[(iv)] If $P$ is generic in $X$ then $f(P \cap X')$ is generic in $Y$;

\item[(v)] If $d=1$ and $f(P \cap X')$ is generic in $Y$ then $P$ is
generic in $X$; and

\item[(vi)] If $f$ respects $P$ and $f(P \cap X')$ is generic in $Y$ then
$P$ is generic in $X$.
\end{itemize}
\end{theorem}
Before proving Theorem~\ref{thm_equivalence}, we emphasise that parts (i)
and (ii) do \textit{not} guarantee that one of the limits involved is
defined exactly if the other is defined. If one of the sequences
converges to some value $c$, then only in the case $c = 0$ can we be
certain that the other will converge. If $c \neq 0$
then the other may fail to converge, although one can easily show that
it will eventually be constrained to vary within the range $[d^{-1} c, d c]$.
We now turn to proving Theorem~\ref{thm_equivalence}.
\begin{proof}
By the bounds on the fibre sizes of $f$ we clearly have 
$$|f(P \cap X' \cap X_n)| \ \leq \ |P \cap X' \cap X_n| \ \leq \ d |f(P \cap X' \cap X_n)|$$
and
$$|f(X' \cap X_n)| \ \leq \ |X' \cap X_n| \ \leq \ d |f(X' \cap X_n)|$$
for all $n \in \mathbb{N}$. It follows from the fact that $f$ is surjective and stratification-preserving
that $f(X' \cap X_n) = Y' \cap Y_n$ and $f(P \cap X' \cap X_n) = f(P \cap X') \cap Y_n$, so
the above inequalities become
$$|f(P \cap X') \cap Y_n| \ \leq \ |P \cap X' \cap X_n| \ \leq \ d |f(P \cap X') \cap Y_n|$$
and
$$|Y' \cap Y_n| \ \leq \ |X' \cap X_n| \ \leq \ d |Y' \cap Y_n|$$
respectively. Now combining these yields
\begin{equation}\label{eq_termineq1}
\frac{1}{d} \frac{|f(P \cap X') \cap Y_n|}{|Y' \cap Y_n|} \ \leq \ 
\frac{|P \cap X_n \cap X'|}{|X_n \cap X'|} \ \leq \ 
d \frac{|f(P \cap X') \cap Y_n|}{|Y' \cap Y_n|}.
\end{equation}
It follows also that
\begin{equation}\label{eq_termineq2}
\frac{1}{d} \frac{|P \cap X_n \cap X'|}{|X_n \cap X'|} \ \leq \ 
\frac{|f(P \cap X') \cap Y_n|}{|Y_n \cap Y'|} \ \leq \ 
d \frac{|P \cap X_n \cap X'|}{|X_n \cap X'|}
\end{equation}
where the left-hand [respectively, right-hand] inequality is obtained by
dividing [multiplying] both sides of the right-hand [left-hand] inequality
in \eqref{eq_termineq1} by $d$.

Now since $X'$ and $Y'$ are generic in $X$ and $Y$ respectively,
Lemma~\ref{lemma_genericsubsetofgeneric} gives
$$\lim_{n \to \infty} \frac{|P \cap X_n|}{|X_n|} \ = \ \lim_{n \to \infty} \frac{|P \cap X' \cap X_n|}{|X_n \cap X'|}$$
and
$$\lim_{n \to \infty} \frac{|f(P \cap X') \cap Y_n|}{|Y_n|} \ = \ \lim_{n \to \infty} \frac{|f(P \cap X') \cap Y_n \cap Y'|}{|Y_n \cap Y'|} = \ \lim_{n \to \infty} \frac{|f(P \cap X') \cap Y_n|}{|Y_n \cap Y'|}$$
where the second equality on the second line holds because $f(P \cap X') \subseteq Y'$.
It is now clear that parts (i) and (ii) follow from 
\eqref{eq_termineq1} and \eqref{eq_termineq2} respectively.

If $f(P \cap X')$ is negligible in $Y$ then the left and right-hand sides of
(i) converge to $0$, from which it follows that the
middle expression converges to $0$, and so $P$ is negligible. Conversely,
if $P$ is
negligible then exactly the same argument applies with (ii) in place of (i)
to show that $f(P \cap X')$ is negligible.
This proves part (iii).

If $P$ is generic in $X$ then $X \setminus P$ is negligible in $X$, so by
part (iii), $f((X \setminus P) \cap X')$ is negligible in $Y$.
But by surjectivity, we must have
$$Y' \setminus f(P \cap X') \subseteq f((X \setminus P) \cap X')$$
so that $Y' \setminus f(P \cap X')$ is negligible in $Y$. Since $Y'$ is
generic in $Y$ and generic sets are closed under intersection, it follows that
$$Y \setminus f(P \cap X') = (Y' \setminus f(P \cap X')) \cup (Y \setminus Y')$$
is negligible in $Y$, so that $f(P \cap X')$ is generic in $Y$ as required
to prove part (iv).

If $d = 1$ and $f(P \cap X')$ is generic in $Y$ then
it is immediate from part (i) that $P$ is generic in $X$, so
that part (v) holds.

Finally, suppose that $f$ respects $P$ and that $f(P \cap X')$ is generic
in $Y$. Since $f$ is surjective we have
$$Y' \ = \ f(X') \ = \ f((X \setminus P) \cap X') \cup f(P \cap X').$$
Now since $f$ respects $P$, we know that $f((X \setminus P) \cap X')$ and
$f(P \cap X')$ are disjoint, and since $Y'$ is generic in $Y$ is follows
that
$$f((X \setminus P) \cap X') \ = \ Y' \setminus f(P \cap X')$$
is negligible in $Y$. But now by part (iii), we deduce that $X \setminus P$
is negligible in $X$, and hence that $P$ is generic in $X$, as required
to prove part (vi).
\end{proof}

A particularly useful special case is the following immediate corollary.
\begin{corollary}\label{cor_equivalence}
Let $X$ and $Y$ be stratified sets, $X' \subseteq X$ and $Y' \subseteq Y$
be generic subsets of $X$ and $Y$ respectively, $f : X' \to Y'$ a surjective,
stratification-preserving, bounded-to-one map. Then for any $P \subseteq X$
such that $f$ respects $P$, we have that $P$ is generic [respectively, negligible]
in $X$  if and only if $f(P \cap X')$ is generic [negligible] in $Y$.
\end{corollary}

Next, we apply Theorem~\ref{thm_equivalence} to show that the generic properties
of finitely presented semigroups are essentially governed by those of
finitely presented monoids. Recall that if $S$ is a semigroup then $S^1$
denotes the monoid with set of elements $S \cup \lbrace 1 \rbrace$ where $1$
is a new symbol not in $S$, and multiplication defined by
$$st = \begin{cases} 
  \text{the $S$-product } st &\text{ if } s, t \in S; \\
  s                  &\text{ if } t = 1; \\
  t                  &\text{ if } s = 1.
\end{cases}
$$

\begin{theorem}\label{thm_monoidvssemigroup}
Let $\mathscr{C}$ be a class of monoids, $A$ a finite alphabet and $k \in \mathbb{N}$.
Then the generic $A$-generated $k$-relation monoid (with respect to either the sum length
stratification or the maximum
length stratification) belongs to $\mathscr{C}$ if and only if the
generic $A$-generated $k$-relation semigroup $S$ (with respect to the
corresponding stratification) is such that $S^1$ belongs to $\mathscr{C}$.
\end{theorem}
\begin{proof}
Let $X$ and $Y$ be the sets of
$k$-relation monoid and semigroup presentations respectively over $A$.
Suppose $X$ and $Y$ are equipped with either the sum length or
the maximum length stratification. Let $P$ be the set of presentations
in $X$ such that the monoid presented lies in $\mathscr{C}$, and let $Q$
be the set of presentations in $Y$ such that the semigroup $S$ presented
is such that $S^1$ lies in $\mathscr{C}$.

Let $Y' = Y$ and let $X' = Y \subseteq X$ be the set of semigroup
presentations viewed as a subset of the set of monoid presentations,
that is, those monoid presentations in which no relation word is
empty.
By Lemma~\ref{lemma_fewshortwords_sumlength} (for the sum length stratification) or
Lemma~\ref{lemma_fewshortwords_maxlength} (for the maximum length stratification)
$X'$ is generic in $X$,
and obviously $Y' = Y$ is generic in $Y$.

Define $f : X' = Y \to Y' = Y$ to be the identity function.
Then $f$ is $1:1$, surjective onto $Y'$, and preserves the sum length and
maximum length stratifications. Letting $d = 1$ and $k_n = 1$ for all $n$,
we see that the conditions of Theorem~\ref{thm_equivalence} are satisfied,
so $P$ is generic in $X$ if and only if $f(P \cap X')$ is generic
in $Y$.

Since $f$ is the identity function on $X'$, a semigroup presentation
$\mathscr{P} \in f(P \cap X')$ 
exactly if $\mathscr{P}$ interpreted as a monoid presentation lies in
$P$. Since $\mathscr{P}$ has no empty relation words, it is easy to see that
the monoid presented by $\mathscr{P}$ is isomorphic to $S^1$, where $S$ is
the semigroup presented by $\mathscr{P}$. Thus, $\mathscr{P} \in f(P \cap X')$
if and only if $S^1 \in \mathscr{C}$, that is, if and only if
$\mathscr{P} \in Q$. Hence, $f(P \cap X') = Q$, and so $P$ is generic in $X$
if and only if $Q$ is generic in $Y$, as required.
\end{proof}

\begin{corollary}
For every $m \geq 1$, $k \in \mathbb{N}$ and alphabet $A$ of size at least
$2$, the generic $A$-generated $k$-relation semigroup
(with respect to either the sum length stratification or the maximum
length stratification) satisfies the small overlap condition $C(m)$.
\end{corollary}

An \textit{unordered monoid presentation} consists of a set $A$ of generators
and an (unordered) \textit{set} $R$ of relations, each of which is an
\textit{unordered} pair of words from $A^*$. Equivalence of words is defined
exactly as for ordered presentations (see Section~\ref{sec_monoids}), as are the
sum length and maximum length stratifications on the sets
of $A$-generated presentations with some fixed number $k$ of relations.
There is an obvious map from the ordered to the unordered presentations
over a given alphabet $A$, which simply ``forgets'' the ordering of the
relations and the ordering of the pair of words in each relation, and
discards any ``duplicate'' relations. Unordered semigroup presentations
can of course be defined analogously.

\begin{theorem}\label{thm_ordervsunorder}
Let $\mathscr{C}$ be a class of monoids, $A$ an alphabet and $k$ a
non-negative integer. Then
the generic [negligible] $A$-generated $k$-relation monoid (with
respect to either the sum length stratification or the maximum length stratification)
belongs to $\mathscr{C}$ if and only if the generic [respectively negligible]
$a$-generator $k$-relation unordered monoid (with respect to the corresponding
stratification) belongs to $\mathscr{C}$. The corresponding statement for
semigroups also holds.
\end{theorem}
\begin{proof}
We prove the result for monoids; that for semigroups can be proved in
exactly the same way. Let $X$ be the set
of ordered $k$-relation monoid presentations over $A$, and $Y$ the set of
unordered $k$-relation monoid presentations over $A$. Let $P \subseteq X$ and
$Q \subseteq Y$ be the sets of presentations in $X$ and $Y$ respectively
such that the monoid presentated belongs to $\mathscr{C}$.

Let $X' \subseteq X$
be the set of ordered presentations which do not feature the same relation
twice, or two relations of the form $(u, v)$ and $(v, u)$ for some distinct
words $u$ and $v$. We have
seen that $C(2)$ presentations do not feature the same relation word twice, so $X'$ certainly contains all the $C(2)$
presentations. It follows by Theorem~\ref{thm_generic_sumlength} (for the sum relation
length stratification) or Theorem~\ref{thm_generic_maxlength} (for the maximum relation
length stratification) that $X'$ is generic in $X$. Let $Y' = Y$; then
certainly $Y'$ is generic in $Y$.

Define $f : X' \to Y' = Y$ to be the restriction to $X'$ of the obvious map
described above from ordered to unordered presentations. It is clear
from the definition of $X'$ that $f$ preserves the number of relations in
the presentation and so really does define a map to $Y$, and moreover that
this map is surjective. Since $f$ takes each ordered presentation to an
unordered presentation of the same monoid, it is also obvious that $f$
respects $P$ and maps $P \cap X'$ onto $Q$. It is easily seen
that $f$ preserves both the sum length and the maximum length
stratifications. Moreover, $f$ clearly has fibre size bounded above by
$k! 2^k$. It follows that the conditions of Corollary~\ref{cor_equivalence}
are satisfied, so that $P$ is generic in $X$ if and only if $f(P) = Q$
is generic in $Y$.
\end{proof}
We thus allow ourselves to speak of a generic monoid or generic semigroup, without
worrying about whether the presentation is defined to have a set or a
sequence of relations.

\section{Properties of Generic Monoids and Semigroups}\label{sec_consequences}

In this section we explore some of the consequences of our results for
generic monoids and
semigroups. Recall that a monoid or semigroup is called \textit{$\mathscr{J}$-trivial}
if distinct elements always generate distinct principal ideals. 

\begin{proposition}\label{prop_jtrivial}
Any $C(3)$ semigroup or monoid is torsion-free and $\mathscr{J}$-trivial.
\end{proposition}
\begin{proof}
Let $S$ be a semigroup or monoid with a $C(3)$ presentation $\langle A \mid R \rangle$.
By a result of Remmers \cite{Remmers71}, only finitely many words over the
alphabet $A$ represent the same element of $S$.

Suppose first that $S$ it is not $\mathscr{J}$-trivial, and choose
$a, b \in S$ be distinct elements generating the same ideal. Then in
particular, $a$ is in the ideal generated by $b$, so we have $a = pbq$
for some $p,q \in S$. But also $b$ is in the ideal generated by $a$, so that
 and $b = ras = rpbqs$ for some $r,s \in S$. Now choose words
$\hat{b}, \hat{p}, \hat{q}, \hat{r}, \hat{s} \in A^*$ representing
$b, p, q, r, s \in S$ respectively. Certainly at least one of $\hat{r}$ and
$\hat{s}$ is non-empty, since otherwise we would have $r = s = 1$ so that 
$b = ras = a$. But now it is easily seen that
$(\hat{r}\hat{p})^i \hat{b} (\hat{q}\hat{s})^i$ represents $b$ for every
$i > 0$, contradicting Remmers' result.

Similarly, suppose $a \in S$ is non-identity torsion element. Then there is a
non-empty word $\hat{a} \in A$ representing $a$. But now it is easy to see that
infinitely many powrs of $\hat{a}$ must represent the same element, again
contradicting Remmers' result.
\end{proof}

Combining with our theorem with have the following.
\begin{theorem}
Let $A$ be an alphabet of size at least $2$ and let $k$ be a positive
integer. Then the monoid defined by the generic $A$-generated $k$-relation
presentation (with respect to either the sum length stratification or the maximum length
stratification) is non-trivial, torsion-free and $\mathscr{J}$-trivial. In
particular, it is not a group, an inverse monoid or a regular monoid.
The corresponding statements for semigroups also hold.
\end{theorem}
\begin{proof}
By Theorem~\ref{thm_generic_sumlength} (respectively Theorem~\ref{thm_generic_maxlength} for the other
stratification) the generic $A$-generated $k$-relation presentation satisfies
$C(3)$, and so by Proposition~\ref{prop_jtrivial} the semigroup presented is
torsion-free and $\mathscr{J}$-trivial. If it were trivial then 
every word over the alphabet would have to represent the identity,
contradicting once more Remmers' result mentioned in the proof of the previous
proposition.
\end{proof}

By a recent result of the author, the uniform word problem for $C(4)$ semigroups
is solvable in time linear in the word lengths and polynomial in the
presentation size \cite[Theorem~2]{K_smallover}. Hence, we obtain
\begin{theorem}
Let $A$ be an alphabet of size at least $2$ and let $k$ be a positive integer.
Then the generic $A$-generated $k$-relation
presentation (with respect to either the sum length stratification or the maximum length
stratification) has word problem solvable in linear time. The corresponding
statement for semigroups also holds.
\end{theorem}

Since there is also an algorithm to decide, in (worst-case) polynomial
time whether a given presentation satisfies the condition $C(4)$
\cite[Corollary~5]{K_smallover}, we also obtain
\begin{theorem}
Let $A$ be an alphabet of size at least $2$ and $k$ be a positive integer.
Then the uniform word problem for $A$-generated, $k$-relation monoid
presentations is generically solvable in polynomial time. The corresponding
statement for semigroups also holds.
\end{theorem}

Further work of the author \cite{K_smallover2} has established a number of
automata-theoretic properties of monoids which admit finite presentations
satisfying the condition $C(4)$. It follows from Theorems~\ref{thm_generic_sumlength} and \ref{thm_generic_maxlength}
that the ``generic'' monoid and semigroup will enjoy all these properties.
The following theorem summarises these properties; for brevity we omit
definitions of terms; which can be found in \cite{K_smallover2}.

\begin{theorem}
Let $A$ be an alphabet of size at least $2$ and let $k$ be a positive
integer. Then the monoid defined by the generic $A$-generated $k$-relation
presentation (with respect to either the sum length stratification or the maximum length
stratification) is rational in the sense of \cite{Sakarovitch87}, asynchronous
automatic and word hyperbolic in the sense of \cite{Duncan04}. It also
satisfies an analogue of Kleene's theorem and has a boolean algebra of
rational subsets and decidable rational subset membership problem.
\end{theorem}

\section*{Acknowledgements}

This research was supported by an RCUK Academic Fellowship.
The author would like to thank A.~V.~Borovik and V.~N.~Remeslennikov
for their many suggestions; he also thanks 
the organisers and participants of the \textit{AIM Workshop on
Generic Complexity}, held in Palo Alto in August 2007, where he had
many helpful conversations, and the American
Institute of Mathematics for funding his attendance there.

\bibliographystyle{plain}

\def\cprime{$'$} \def\cprime{$'$} \def\cprime{$'$}

\end{document}